\documentclass{amsart}

\newtheorem{thm}{Theorem}[section]
\newtheorem{proposition}[thm]{Proposition}

\newtheorem{lemma}[thm]{Lemma}
\newtheorem*{thmintro}{Theorem}

\theoremstyle{definition}
\newtheorem{definition}[thm]{Definition}
\newtheorem{example}[thm]{Example}

\theoremstyle{remark}
\newtheorem{remark}[thm]{Remark}

\newcommand{\RR}{\mathbf{R}}
\newcommand{\ZZ}{\mathbf{Z}}
\newcommand{\NN}{\mathbf{N}}

\newcommand{\F}{\ensuremath{\mathcal{F}}}

\newcommand{\CH}{\ensuremath{\mathcal{H}}}
\newcommand{\CD}{\ensuremath{\mathcal{D}}}

\newcommand{\eps}{\varepsilon}

\newcommand{\hess}{\mathrm{Hess}}
\newcommand{\Fix}{\mathrm{Fix}}
\begin{document}

\title{Singular Riemannian Foliations on Nonpositively Curved Manifolds}





\title{Singular Riemannian Foliations on Nonpositively Curved Manifolds}
\author{Dirk T\"oben}
\address{Mathematisches Institut, Universit\"at zu K\"oln, Weyertal 86-90, 50931 K\"oln, Germany}
\email{dtoeben@math.uni-koeln.de}
\date{September 12, 2005}
\begin{abstract}
We prove the nonexistence of a proper singular Riemannian foliation admitting section in compact manifolds of nonpositive curvature. Then we give a global description of proper singular Riemannian foliations admitting sections on Hadamard manifolds. In addition by using the theory of taut immersions we provide a short proof of this result in the special case of a polar action. 
\end{abstract}

\maketitle
\section{Introduction}
A singular Riemannian foliation with sections (see definition \ref{srf}) can be seen as a generalization of an isoparametric family in a Euclidean space to arbitrary ambient spaces. Another important example is the orbit decomposition of a polar action (see definition \ref{dfnpolar}). These actions are well understood in compact symmetric spaces; see \cite{kollross} for a classification of hyperpolar actions and \cite{podestathorbergsson} for a classification of polar actions in compact symmetric spaces of rank one. (For surveys on these topics, see \cite{thorbergsson1} and \cite{thorbergsson2})\par
No classification of polar actions is known for symmetric spaces of noncompact type. Before we state our result, we give some examples. If $(G,K)$ is a symmetric pair of noncompact type, then the left action of $K$ on $G/K$ is hyperpolar as in the compact case and has compact, isoparametric orbits. If $G=KAN$ denotes the Iwasawa decomposition, then the left action of $N$ on $G/K$ is hyperpolar with noncompact and regular orbits. Now let $M$ be a closed submanifold of the hyperbolic space $H^n$ such that its principal curvatures are not larger than $1$. Then the set of distance tubes around $M$ defines a singular Riemannian foliation $\F_M$ with normal geodesics of $M$ as sections. One can easily choose $M$ such that $\F_M$ is inhomogeneous, i.e. does not come from a polar action. This example can be easily generalized to Hadamard manifolds. The following theorem gives a global description of a proper singular Riemannian foliation admitting sections in an Hadamard manifold. Please note how the above examples are covered.
\begin{thmintro}
Let $\F$ be a proper singular Riemannian foliation with sections in an Hadamard manifold $X$. Then $\F$ is the product foliation of a compact isoparametric foliation and the trivial foliation of $\RR^n$ with $\RR^n$ as the only leaf.
\end{thmintro}
This is Theorem \ref{mainthmhadamard}, proved in section \ref{sec:srf}. The homogeneous counterpart is the following:
\begin{thmintro}
Let $G$ be a Lie group acting polarly and properly on an Ha\-da\-mard manifold $X$ and $K$ be a maximal compact subgroup of $G$. Then there is a natural $G$-equivariant diffeomorphism $X\cong G\times_K \RR^s, s =\dim X-\dim G/K$, where the action of $K$ on $\RR^s$ is orbit equivalent to an $s$-representation.
\end{thmintro}
This is Theorem \ref{mainthm1}. We give a short proof of it in section \ref{sec:polar}. Now what happens if the ambient space is a compact quotient of an Hadamard manifold?
\begin{thmintro}
Except possibly for a regular Riemannian foliations there are no proper singular Riemannian foliations admitting sections on a compact, nonpositively curved Riemannian manifold $N$. 
\end{thmintro}
This is Theorem \ref{mainthmcpt} proved in section \ref{sec:srf}.\par
I would like to thank Professor Thorbergsson for many helpful discussions and the Deutsche Forschungsgemeinschaft for their support by the Schwerpunktsprogramm "Globale Differentialgeometrie" (SPP 1154).

\section{Polar Actions}\label{sec:polar}
\begin{definition}\label{dfnpolar}
An isometric action of a Lie group $G$ on a Riemannian manifold is {\it polar}, if there is an immersed submanifold $\Sigma$, called {\it section}, that meets any orbit and is orthogonal to them at each point of intersection.
\end{definition}
Let $G$ be a Lie group acting by isometries on a complete Riemannian manifold $N$. The action is {\it variationally complete}, if for any orbit $M$ and any normal geodesic $\gamma$, an $M$-Jacobi field along $\gamma$ that is tangent to some other orbit is the restriction of a $G$-Killing field to $\gamma$. Now let $G$ be a Lie group acting properly and polarly on an Hadamard manifold $X$. For this case the following theorem applies.
\begin{thm}[Conlon]
A proper polar action whose sections have no conjugate points is variationally complete.
\end{thm} 
The proof of the theorem is essentially due to Conlon \cite{conlon}. He proves it for flat sections, but it also works for (not necessarily embedded) sections without conjugate points. Also the compact action group is replaced by a proper action.\par
Let $M$ be a closed and embedded submanifold in $X$ and let $q\in X$ be a point that is not focal point of $M$. Then the squared distance funtion $f_q:=d^2(\cdot,q):M\to \RR$ is a proper Morse function. Let $M_s=M\cap f_q^{-1}([0,s))$ for any $s\in [0,\infty)$. The Morse inequalities say $b_i(M^s)\leq \mu_i(f_q|M_s)$, where $b_i(M^s)$ is the $i$-th Betti number of $M^s$, say over $\ZZ_2$, and $\mu_i(f_q|M^s)$ is the number of critical point of index $i$ of $f_q|M^s$. We say that $M$ is {\it taut} in $X$, if all the Morse inequalities are equalities for any nonfocal point $q\in X$ and any $s$ (note that we replaced the energy functional $E_q$ defined on an infinite dimensional path space by $f_q$). By \cite{bottsamelson} a variationally complete action has taut orbits:
\begin{thm}[Bott-Samelson]
Each orbit $M$ of $G$ is taut in $X$.
\end{thm}
Using this result we can give a geometric description of a proper polar actions of a Lie group $G$ on an Hadamard manifold $X$. Let $K$ be a maximal compact subgroup of $G$. By Cartan's Theorem $K$ fixes some point $p\in X$, so $K\subset G_{p}$ and by maximality $K=G_{p}$. It is known that $M:=Gp\cong G/K$ is diffeomorphic to $\RR^m$, where $m=\dim Gp$. Therefore $b_i(Gp)$ is $1$ for $i=0$, otherwise $0$. As $M$ is taut, any nonfocal point of $M$ in $X$ has exactly one preimage under $\eta:=\exp^\perp:\nu M\to X$. This implies that there are no focal points along a normal geodesic of $M$. So $M$ has no focal points at all. Thus there is exactly one normal geodesic segment ending in a given point $q\in X$. Therefore the normal exponential map $\eta$ is a diffeomorphism. Let $S:=\exp(\nu_pGp)$ and $s:=\dim S$. Now $X$ is a tube of $Gp$ with section $S$. By the Tube Lemma following two maps
$$
G\times_K \nu_{p}Gp\longrightarrow \nu Gp\stackrel{\exp^\perp}{\longrightarrow} X,
$$
are G-equivariant diffeomorphisms. The action of $K$ on $\RR^s$ is the linearization of the $K$-action on $S$; it is the polar slice representation of $G$ at $p$ and therefore by Dadok (\cite{dadok}) orbit equivalent to an $s$-representation, i.e. the isotropy representation of a symmetric space. We sum up:
\begin{thm}\label{mainthm1}
Under the above assumptions, there is a natural $G$-equivariant diffeomorphism $X\cong G\times_K \RR^s$. The action of $K$ on $\RR^s$ is orbit equivalent to an $s$-representation.
\end{thm}
Abels has proven a topological version of the first statement in open manifolds, see Corollary 1.3 in \cite{abels}, in which $\RR^s$ is replaced by a slice for the $G$-action. He obtains his result by constructing this slice. 
It is natural to ask whether the slice is canonical, i.e. of the form $\exp(\nu_{p}Gp),$ in the case of an isometric action. This is not true as the one parameter group of glide rotations along an axis in $\RR^3$ shows. But for polar actions these canonical slices exist. This is the main point of this theorem.\par
The theorem tells that by projection of $X=G\times_K \RR^s$ onto the first factor $G/K=Gp$ every orbit $G(\exp v), v\in \nu_pGp$ fibers over $Gp$ with fiber $Kv$ and structure group $K$. This fiber bundle is trivial, since $G/K\cong \RR^s$, so $G(\exp v)\cong \RR^s\times Kv$; this is a cylinder over an orbit of an s-representation, which is an isoparametric submanifold. We remark that is it not difficult to show that this last cylindric decomposition is metric if $X$ is the Euclidean space.

\section{Singular Riemannian Foliations admitting Sections}\label{sec:srf}
We want to prove an analogous result for singular Riemannian foliations admitting sections. 
\begin{definition}[\cite{molino}]\label{srf}
Let $\F$ be a partition of injectively immersed submanifolds (the {\it leaves}) of a Riemannian manifold $N$. For any $p\in N$ let $L_p$ be the leaf through $p$ and let $T\F=\bigcup_{p\in N}T_pL_p$. We define $\Xi(\F)$ as the module of (differentiable) vector fields on $N$ with values in $T\F$. We call $\F$ a {\it singular foliation}, if it satisfies the following differentiability condition, a {\it singular Riemannian foliation}, if it satisfies in addition the transnormality condition:
\begin{enumerate}
\item (Differentiability) $\Xi(\F)$ acts transitively on $T\F$, i.e., for any $v\in T_p\F, p\in N$ there is $X\in\Xi(\F)$ with $X_p=v$.
\item (Transnormality) a geodesic starting orthogonally to a leaf intersects the leaves it meets orthogonally;
\end{enumerate}
A leaf of maximal dimension is called {\it regular}, and so each point of it, otherwise {\it singular}. If, in addition, for any regular $p$ there is an isometrically immersed complete totally geodesic submanifold $\Sigma_p$ (the {\it section}) with $T_p\Sigma=\nu_pL_p$, that meets any leaf and always orthogonally, $\F$ is a singular Riemannian foliation {\it admitting sections}.
\end{definition}
We call a singular Riemannian foliation {\it proper} if all leaves are properly immersed, i.e. closed and embedded. 
\begin{example} The orbit decomposition of a polar action is a singular Riemannian foliation admitting sections. We give another important example. Let $M$ be an isoparametric submanifold in $\RR^n$. Then the partition of parallel and focal submanifolds of $M$ is a singular Riemannian foliation with sections. We call this partition an {\it isoparametric foliation}. This special kind of singular foliation occurs in every singular Riemannian foliation with sections as we will see in the next theorem. We remark that isoparametric submanifolds in $\RR^n$ have been generalized to {\it equifocal submanifolds} in simply connected compact symmetric spaces in \cite{terngthorbergsson}, to {\it submanifolds with parallel focal structure} in general symmetric spaces in \cite{ewert} and to arbitrary ambient spaces in \cite{alexandrino} and \cite{toeben}.
\end{example}
Now let $\F$ be a singular Riemannian foliation with sections in a complete Riemannian manifold $N$. Let $x\in N$ and $P$ be an open, relatively compact neighborhood of $x$ in $L_x$. Then there is an $\eps>0$ such the restriction of $\exp^\perp:\nu L_x\to N$ to $\{X\in\nu P\mid \|X\|<\eps\}$ is a diffeomorphism onto its image $T$. Let $\pi:T\to P$ the be orthogonal projection. Condition (1) of a singular Riemannian foliation implies that by eventually shrinking $\eps$ we can assume that the leaves intersect the {\it slices} $S_y:=\pi^{-1}(y),y\in P$ transversally. Then $T$ is a {\it distinguished neighborhood} in the sense of Molino (\cite{molino}). The transversal intersections induce a singular foliation on a slice $S_x$, such that each leaf has the same codimension as its restriction to $S_x$. We can lift this singular foliation to a singular foliation on a ball neighborhood in $TS_x=\nu_xL_x$ of the origin via the exponential map. This singular foliation is homothety invariant by the Homothety Lemma (Lemma 6.2 in \cite{molino}) applied to $\F$ at $x$; therefore we can extend it to $\nu_xL_x$. The singular foliation obtained is denoted by $\F_x$. The following theorem is due to Alexandrino (Theorem 2.10 in \cite{alexandrino}):
\begin{thm}[Slice Theorem]
$\F_x$ is an isoparametric foliation. The sections are exactly of the form $T_x\Sigma$, where $\Sigma$ is a section of $\F$ through $x$.
\end{thm}
Now consider a section $T_x\Sigma$ for $\F_x$. The set of singular points in $T_x\Sigma$ is a finite union of hyperplanes through the origin. The {\it generalized Weyl group} $W$ is generated by the reflections across these hyperplanes.\par 
Let $\F$ be a singular Riemannian foliation admitting sections of an Ha\-da\-mard manifold $X$. A section $\Sigma$ is also an Ha\-da\-mard manifold. It is moreover closed and embedded in $X$. We want to introduce an analogue of the generalized Weyl group of a polar action. For details see \cite{toeben}. Let $M$ be a regular leaf. Now let $\tau$ be an arbitrary curve in $M$ starting and ending in $M\cap\Sigma$. We obtain a map $T_{\tau(0)}\Sigma=\nu_{\tau(0)}M\to \nu_{\tau(1)}M=T_{\tau(1)}\Sigma$. By exponentiating we obtain an isometry from a ball neighborhood of $\tau(0)$ in $\Sigma$ to a ball neighborhood of $\tau(1)$ in $\Sigma$. This map preserves leaves. One can extend this map to an isometry $\phi_\tau:\Sigma\to\Sigma$ preserving leaves. The collection of these $\phi_\tau$ over all such curves $\tau$ is a group, called {\it transversal holonomy group}. 
$$
\Gamma\subset I(\Sigma)\quad\mbox{and}\quad\Gamma p=L_p\cap \Sigma\ \mbox{for all}\ p\in\Sigma.
$$
This group is independent of the choice of $M$. The linearization $d\Gamma_q$ of the isotropy group $\Gamma_q$ of a regular point $q\in\Sigma$ is the normal holonomy of $L_q$ and the orbit of any $q\in\Sigma$ describes the recurrence of $L_q$ to the section $\Sigma$. This property characterizes $\Gamma$ but also the generalized Weyl group of a polar action. This implies the equality of both notions for polar actions.\par
Now let $\F$ be proper, so any leaf is closed and embedded. This implies that $\Gamma$ is discrete. Let $M$ be a regular leaf with trivial normal holonomy. Let $\{p_i\}_{i\in I}:=M\cap \Sigma$. We define the set $\{\CD_{p_i}\}_{i\in I}$ of (closed) Dirichlet domains by
$$
{\CD}_{p_i}=\{q\in \Sigma\ |\ d(p_i,q)\leq d(p_j,q)\ \mbox{for all}\ j\neq i\}.
$$
The group $\Gamma$ acts transitively on the set $\{\CD_{p_i}\}_{i\in I}$. Since $M$ has trivial holonomy we have $\Gamma_{p_i}=\{e\}$, and $\Gamma$ acts simply transitively on this set. \par
We take a look at the boundary of a Dirichlet domain. For $x,y\in\Sigma$ let $H_{x,y}:=\{z\in\Sigma\mid d(x,z)=d(y,z)\}$ be the bisector/central hypersurface between $x$ and $y$. In an Hadamard it is a hypersurface that disects $\Sigma$. The set $F=H_{p_i,p_j}\cap \CD_{p_i}$ is called a {\it wall} of $\CD_{p_i}$, if it contains an open non-empty subset of $H_{p_i,p_j}$. Two Dirichlet domains are called {\it neighbors}, if they contain a common wall.\par
Now we fix one $p:=p_i$ and let $\CD=\CD_{p_i}$. 
\begin{proposition}
Let $\F$ be a proper singular Riemannian foliation admitting sections on an Hadamard manifold $X$. Let $\Sigma$ be a section of $\F$ and $\Gamma$ be the transversal holonomy group acting on $\Sigma$. Then $\Gamma$ is a Coxeter group.
\end{proposition}
\begin{proof}
By the Poincar\'e-Lemma (see for instance Lemma 2.5 in \cite{aklm}), $\Gamma$ is generated by the set of elements of $\Gamma$ that map $\CD$ to a neighboring domain. We want to show that these elements are reflections. This implies that $\Gamma$ is a Coxeter group by Theorem 3.5. in \cite{aklm} (see the definition of a Riemannian Coxeter manifold in subsection 3.2. of \cite{aklm}). Let $\CD'$ be a neighboring Dirichlet domain of $\CD$ and $g\in\Gamma$ be the unique element with $ g(\CD)=\CD'$. Then $\CD$ and $\CD'$ have a common wall $F$. We claim that the wall $F$ consists of singular points of $\F$. If not, there is a regular point $q\in F$. In \cite{toeben} the following was shown: singular points are focal points for any regular leaf, in particular for $M$; conversely, a focal point $x$, say in the section $\Sigma$, that is not conjugate in $\Sigma$ to any point in $M\cap\Sigma$, is a singular point of $\F$. Since in our case the sections have no conjugate points, the focal points of $M$ are exactly the singular points of $\F$. Thus $q$ is not a focal point and the squared distance function $f_q=d^2(\cdot,q):L_p\to\RR$ is a proper Morse function. As can be read for instance in \cite{toeben}, $f_q$ has only one local minimum (regular leaves are $0$-taut). We have two minimal normal geodesics $\gamma_1$ and $\gamma_2$ of $M$ from $p$ respectively $gp$ to $q$. Thus $p$ and $gp$ are critical points of $f_q$. The segment $\gamma_i$ does not meet any singular points, because the singular points are in the boundary of the Dirichlet domains as shown in \cite{toeben}. In particular there are no focal points of $M$ on $\gamma_i$. Hence $p$ and $gp$ are critical points of $f_q$ of index $0$, i.e. minima, contradiction. So $F$ only consists of singular points. Note that by using $0$-tautness similarly as above one can show that $(X,\F)$ has no exceptional leaves, i.e. regular leaves with nontrivial holonomy (see Lemma 1A.3 in \cite{podestathorbergsson} and \cite{toeben}).\par
We choose $q$ from the interior of $F$. Since $q$ is singular, there has to be a singular hyperplane in $T_q\Sigma$ for the isoparametric foliation $\F_q$. We denote the reflection across this hyperplane by $s$. By definition this element is contained in the Weyl group of $\F_q$ acting on $T_q\Sigma$. As $q$ is in the interior of a wall, there cannot be other singular hyperplanes through the origin of $\nu_qL_q$. Thus $W$ is generated by $s$. Any element of $W$ can be extended to a leaf preserving isometry of $\Sigma$ in $\Gamma$. In this sense we write $W\subset\Gamma$. Then $H_{p,gp}=\Fix (<s>)$ and $g=s$. This shows that $\Gamma$ is a Coxeter group. 
\end{proof}
The proof shows that the bisector $H_{p,gp}$ as the fixed point set of a disecting reflection is a connected and totally geodesic hypersurface. Moreover it only consists of singular points, otherwise they would belong to exceptional leaves). Let $\{A_i\}_{i\in J}$ be the set of reflection hyperplanes of $\Gamma$. Then the union of the $A_i$'s is the set of singular points in $\Sigma$. We call $A_i$ a {\it singular hyperplane}. Each $A_i$ as a bisector of a disecting reflection is totally geodesic and complete in the Hadamard manifold $\Sigma$, therefore closed and embedded (and itself an Hadamard manifold). The connected components of $\Sigma\backslash\bigcup_{i\in I}A_i$ are by definition the chambers of $\Gamma$. By Corollary 3.8 of \cite{aklm} the chambers coincide with the Dirichlet domains of regular $\Gamma$-orbits.\par
The Coxeter groups have been classified for the Euclidean space but not for the hyperbolic space. A priori we can expect complicated Coxeter groups for the transversal holonomy group. But not every Coxeter group can occur. 
\begin{proposition}
Let $\F$ be a proper singular Riemannian foliation admitting sections on an Hadamard manifold $X$. Let $\Sigma$ be a section of $\F$ and $\Gamma$ be the transversal holonomy group acting on $\Sigma$. Then $\Gamma$ is isomorphic to a (finite) Coxeter group of Euclidean space.
\end{proposition}
\begin{proof} We want to show that $|\Gamma|$ is finite. Then it follows from Cartan's Theorem that $\Gamma$ has a fixed point $p_0$. Thus $(T_{p_0}\Sigma,d\Gamma)$ is a finite Coxeter group of Euclidean space and completely describes $(\Sigma,\Gamma)$.\par
Therefore we assume $|\Gamma|=\infty$. Let $C_0=\CD$ be a chamber. Let $S$ be the set of generators of $\Gamma$ given by the reflections across the walls of $C$. We claim that the set $R=\{l(g)\mid g\in\Gamma\}\subset\NN$ is not bounded, where $l(g)$ denotes the length of $g$ with respect to $S$.\par
{\sc Case 1}: $|S|<\infty$. If $R$ were bounded, say by a number $N$, then $|\Gamma|\leq |S|^N<\infty$, contradiction.

{\sc Case 2}: $|S|=\infty$. Let us give some definitions first. A sequence of chambers $C_0,\ldots, C_l$, for which any two consecutive chambers have a common wall, is called a {\it gallery of length} $l$. It is called {\it minimal} if there is no gallery from $C_0$ to $C_l$ of length smaller than $l$. For any chamber $C$, there is a unique $g_C\in \Gamma$ with $g_C(C_0)=C$. This correspondence between chambers and elements of $\Gamma$ is bijective. The length of a minimal gallery from $C_0$ to a chamber $C$ is equal to $l(g_C)$ (see \cite{gutkin} for the simple relation between a gallery and a word in $\Gamma$ with respect to $S$). Now let $C_0,\ldots, C_l$ be a minimal gallery. The chamber $C_i$ is given by the reflection of $C_{i-1}$ across the singular hyperplane $A_i$ that contains a common wall. A gallery is minimal if and only if the $A_i, i=1,\ldots l,$ are pairwise different, since another gallery with the same first and the same last chamber would have to pass each $A_i$ at least once (see Corollary 2 in \cite{gutkin}). As the set of walls of a chamber is infinite by assumption we can choose a singular hyperplane $A_{l+1}$ that contains a wall of $C_l$ and is different from the $A_i,i=1,\ldots,l,$. Reflecting $C_l$ across $A_{l+1}$ gives a chamber $C_{l+1}$ which extends the given gallery to a minimal gallery of length $l+1$.
This implies that $R$ is unbounded.\par
So in any case $R$ is unbounded; thus there is a $g\in\Gamma$ of length larger than $n$, the dimension of a regular leaf. Let $\Gamma$ be a geodesic from a regular point in $C_0$ to a point in $gC_0$. By slightly perturbing the endpoints, we can assume that $\gamma$ meets the singular hyperplanes only one at a time. Then $\gamma$ intersects at least $l(g)>n$ singular hyperplanes. This contradicts the following proposition, since $\gamma$ is an $L_{\gamma(0)}$-geodesic.
\end{proof}

\begin{proposition}\label{finitefocal}
Let $M$ be a closed and embedded submanifold of an Hadamard manifold $X$ and $\gamma$ be a normal ray of $M$. Then the number of focal points along $\gamma$, counted with multiplicity is not larger than $\dim M$.
\end{proposition}
\begin{proof}
We assume that $\gamma$ is a unit speed geodesic. Let $x\in M$ be the initial point of $\gamma$. We define $f_t:M\to\RR$ by $f_t(y)=d(y,\gamma(t)).$ The sum of the focal points of $M$ along $\gamma|[0,t]$, counted with multiplicity, is equal to the index of $\hess f_t(x)\subset T_xM$, which is not larger than $\dim M$.
\end{proof}

The following lemma is probably well-known.
\begin{lemma}\label{isoparametric}
Let $\F$ be an isoparametric foliation on $\RR^n$ with compact leaves and $\Sigma$ be a section on which the generalized Weyl group $W$ acts. Then $\Fix(W)$ is equal to the intersection $I$ of all sections through the origin and to the stratum $S$ of zero-dimensional leaves.
\end{lemma}
\begin{proof}
We can assume that the leaves are contained in spheres around the origin. We can decompose $(\RR^n,\F)=(V,\F_1)\oplus (V^\perp,\F_2)$, where $V$ is a subspace of $\RR^n$ and $\F_1$ is an isoparametric foliation, in which one and therefore all regular leaves are full (i.e. not contained in a proper subspace), and where $\F_2$ is the trivial foliation of $V^\perp$ by points.
We want to show the chain of inclusions
$$
V^\perp\subset I\subset S\subset\Fix(W)=V^\perp. 
$$
Let $M$ be a regular leaf of $\F_1$; then $M$ is also a regular leaf of $\F$. The section of $\F_1$ through $q\in M$ has the form $q+\nu^V_qM$, where $\nu^V_qM$ is the normal space of $M$ in $V$ at $q$;  the section of $\F$ through $q\in M$ has the form $q+\nu^{\RR^n}_qM$, where $\nu^{\RR^n}_qM$ is the normal space of $M$ in $\RR^n$ at $q$. We have the direct orthogonal sum $\nu^{\RR^n}_qM=\nu^V_qM\oplus V^\perp$. Since $M$ meets every section, $V^\perp\subset I$. Now $I\subset S$ since the dimension of the family of sections through a point $x\in \RR^n$ considered as a submanifold in $G_k(\RR^n), k=\dim \Sigma,$ gives the {\it defect} $\dim M-\dim L_x$. The zero-dimensional leaves meet $\Sigma$ exactly once. Therefore $S\subset \Fix W$. As $W$ acts trivially on $V^\perp\subset \Sigma$ and $\Fix(W|V\cap\Sigma)=\{0\}$ by Proposition 6.2.3. in \cite{palaisterng}, we have $\Fix(W)=V^\perp$.
\end{proof}
\begin{remark}
The proof shows that $\F_1$ has only one zero-dimensional leaf, name\-ly the origin, which is the intersection of all sections of $\F_1$.
\end{remark}

Now we come back to the case $(X,\F)$.
\begin{proposition}\label{stratum}
Let $S$ be the stratum of leaves of $\F$ with least dimension, say $m$, let $\Sigma$ be a section and $\Gamma$ be the transversal holonomy group acting on $\Sigma$. Then $S\cap\Sigma=\Fix(\Gamma)$ and consequently this is a connected, totally geodesic submanifold. Any leaf $M$ in this stratum is diffeomorphic to $\RR^m$ and $\exp^\perp:\nu M\to X$ is a diffeomorphism.
\end{proposition}
\begin{proof}
Let $p\in\Fix(\Gamma)$ and $x\in S\cap\Sigma$. We know that $\F$ induces an isoparametric foliation $\F_x$ on $\nu_xL_x$. Let $W$ be its generalized Weyl group acting on $T_x\Sigma$ and extend each element of $W$ to an isometry of $\Sigma$. We have $W\subset\Gamma_x\subset \Gamma_p=\Gamma$. Thus the intersection of all sections through $x$, being equal to $\Fix(W)$ by the second statement of the Slice Theorem and Lemma \ref{isoparametric}, also contains the geodesic $\gamma_{px}$. Therefore the family of sections through $p$ has the same dimension as that of $x$, so $p\in S$. We have shown $\Fix(\Gamma)\subset S\cap\Sigma$. We prove $\Fix(\Gamma)\supset S\cap\Sigma$ later. \par
Let $M$ be a leaf through a fixed point $p$ of $\Gamma$. We already know that $M$ has the lowest dimension among all leaves. We want to show that its normal exponential map $\eta:=\exp^\perp:\nu M\to X$ is a diffemorphism. Surjectivity follows from properness. We prove injectivity. Assume $\eta(v_1)=\eta(v_2)=q$ for $v_1,v_2\in \nu M$ with foot points $x_1,x_2$. Let $\Sigma_1$ respectively $\Sigma_2$ be a section to which $v_1$ respectively $v_2$ is tangential. Both sections contain $q$. Let $\Gamma^i$ denote the transversal holonomy group acting on $\Sigma_i$; then $x_i$ is a fixed point of $\Gamma^i$, as two transversal holonomy groups acting on $\Sigma$ respectively $\Sigma'$ are conjugate by a leaf-preserving isometry $\Sigma\to\Sigma'$. Let $W$ be the generalized Weyl group of $\F_q$ acting on $T_q\Sigma_1$. As before, we consider $W\subset \Gamma^1_q$. Now $W\subset\Gamma^1_q\subset\Gamma^1_{x_1}=\Gamma^1$. This means that the geodesic $\gamma_{x_1q}$ from $x_1$ to $q$ is fixed under $\Gamma^1_q$ and therefore also under $W$. But $\Fix(W)$ is the intersection of all sections through $q$. In particular $\gamma_{x_1q}$ is contained in both $\Sigma_1$ and $\Sigma_2$. Thus $x_1\in\Sigma_2$. So $M$ intersects $\Sigma_2$ in $x_1$ and $x_2$. But $M\cap\Sigma_2=\Gamma^2(x_2)=\{x_2\}$, hence $x_1=x_2$. Since $\Sigma_2$ is an Hadamard manifold, $\gamma_{x_1q}=\gamma_{x_2q}$ and therefore $v_1=v_2$ and we have proved injectivity of $\eta$. \par

Injectivity of $\eta$ also implies that $M$ has no focal points: assume $v\in \nu M$ is a singular point for $\eta$, then $\gamma_v|[0,t]$ is not a minimal normal geodesic of $M$ for $t>1$. Fix such a $t$. Thus there is also a minimal normal geodesic of $M$ ending in $\gamma_v(t)$ which contradicts that $\eta$ is injective.\par
This means that $\eta:\nu M\to X$ is a diffeomorphism. We fix a point $q\in X$. Then the squared distance function $f_q=d^2(\cdot,q):M\to\RR$ has only one minimum, say $x$, and no other critical points. By standard Morse theory we can use the flow of the negative gradient field to define a diffeomorphism of $M$ to a ball neighborhood of $x$. Therefore $M^m$ is diffeomorphic to $\RR^m$.\par
Now we want to show $S\cap\Sigma\subset \Fix(\Gamma)$. Let $q\in S\cap\Sigma$. We choose $p\in\Fix(\Gamma)\subset S\cap\Sigma$ and define $M=L_p$. As we have shown the normal exponential map of $M$ is a diffeomorphism. In the words of Molino $X$ is a (global) distinguished neighborhood of $M$. Therefore any leaf meets the slices of $M$ transversally. The restriction $\rho$ of the orthogonal projection $X\to M$ to $L_q$ is a surjective diffeomorpism. We want to show that $\rho$ is a diffeomorphism. For $p'\in M$ we have $\rho^{-1}(p')=S_{p'}\cap L_q$, where $S_{p'}$ is the global slice of $M$ through $p'$. By Proposition 2.1.(a) in \cite{alexandrino} $S_{p'}$ is the union of all sections through $p'$. Let $\Sigma'$ be such a section. Then $L_q\cap\Sigma'=\Gamma'(q')$, where $\Gamma'$ is the transversal holonomy acting on $\Sigma'$ and $q'\in L_q\cap\Sigma'$. If a point of $\Gamma'(q')$ were not in the intersection of all sections through $p'$, then $\dim (S_{p'}\cap L_q)>0$, contradiction. So $S_{p'}\cap L_q$ is in this intersection and is therefore equal to $\Gamma'(q')$. This implies that $\rho:L_q\to M$ is a finite covering with degree $|\Gamma'(q')|=|\Gamma(q)|$. But since $M\cong \RR^m$ the map $\rho$ must be a diffeomorphism. In particular $\Gamma(q)=\{q\}$, so $q\in\Fix(\Gamma)$.  
\end{proof}
Let $N$ be compact Riemannian manifold of negative curvature. Then its isometry group is trivial and in particular $N$ admits no polar actions. In \cite{walczak} it was shown that $N$ admits no Riemannian foliation. Is this also true for singular Riemannian foliations? The theorem below gives a partial answer.
\begin{thm}\label{mainthmcpt}
Except possibly for regular Riemannian foliations there are no proper singular Riemannian foliations admitting sections on a compact, nonpositively curved Riemannian manifold $N$. 
\end{thm}
\begin{proof}
Assume there is a singular Riemannian foliation $\F$ on $N$ that admits sections and has singular leaves. Then we can lift it to $\tilde\F$ to the Hadamard manifold $\tilde N$ along the Riemannian universal covering $\pi:\tilde N\to N$. Clearly $\tilde\F$ is a singular Riemannian foliation admitting sections. Let $M$ be a leaf of $\tilde\F$ of least dimension and $p\in M$. Let $B$ be a ball around $p$ whose radius $r$ is larger than the diameter of $N$; this implies that the translations of $B$ by the action of $\pi_1(N)$ exhausts $\tilde N$. The deck transformation group $\pi_1(N)$ respects $\tilde\F$ by definition. This implies that $M$ is mapped onto a leaf of the singular stratum $S$ of leaves of lowest dimension $m$. Let $\Sigma$ be a section through $p$. By Proposition \ref{stratum} $S\cap \Sigma=\Fix(\Gamma)$ is a connected, totally geodesic submanifold of $\Sigma$. It has at least codimension one because $\F$ has regular {\it and} singular leaves. Thus there is a regular point $q\in\Sigma$ with $d_\Sigma(q,S\cap\Sigma)>r$. The distance from $q$ to the orbit $\pi_1(N)p$ attains its minimum at a point $p'\in\pi_1(N)p$. We have
$$
d(q,\pi_1(N)p)=d(q,p')>d(q,L_{p'}).
$$
There is a point $p''\in L_{p'}$ realizing this distance, i.e. $d(q,p'')=d(q,L_{p'})$. This point $p''$ has to be in $\Sigma$; it is the unique point of the intersection $L_{p'}\cap\Sigma\subset S\cap\Sigma$. Now $d(q,p'')>d_\Sigma(q,S\cap\Sigma)>r$. Altogether $d(q,\pi_1(N)p)>r$, so the action of $\pi_1(N)$ on $B$ does not cover $q$, contradiction.
\end{proof}
Now we prove the analogue of Theorem \ref{mainthm1} for singular Riemannian foliations admitting sections. 
\begin{thm}\label{mainthmhadamard}
Let $\F$ be a proper singular Riemannian foliation with sections in an Hadamard manifold $X$. Then $\F$ is the product foliation of a compact isoparametric foliation and the trivial foliation of $\RR^m$ with $\RR^m$ as the only leaf.
\end{thm}
\begin{proof}
Let $M$ be a leaf of lowest dimension $m$. The normal exponential map $\exp^\perp:\nu M\to X$ is a diffeomorphism. We denote the orthogonal projection onto $M$ by $\rho: X\to M$. Let $\CH$ be the horizontal distribution of $\rho$. We now want to change the metrics on $M$ and $X$ such that $\rho$ becomes a Riemannian submersion. First, since $M\cong\RR^{m}$, we can introduce a flat metric on $M$. Moreover we change the metric on $X$ by keeping its induced metric on the fibres of $\rho$. We demand that the fibres are orthogonal to $\CH$. Now we change the metric on the horizontal distribution $\CH$ of $\rho$ such that $\rho$ becomes a Riemannian submersion. For any point $p\in M$ let $S_p$ be the global slice $\exp(\nu_pM)$ of $\F$. For any $q\in S_p$ let $S_q$ be a local slice. Then $S_q\subset S_p$ (see Proposition 2.1.(b) in \cite{alexandrino}, also \cite{toeben}), so $\CH_q\subset T_qL_q$. This means that horizontal lifts of curves in $M$ remain in leaves of $\F$. We fix $p\in M$.
We identify $M$ and $T_pM$ and we define a map $T_pM\times \nu_pM\rightarrow X$ by mapping $(v,w)$ to the endpoint of the lift of the straight line $\gamma_w$ in $T_pM$ along $\rho$ to the starting point $\exp^{S_p}(v)\in S_p$. This map $\phi$ is a diffeomorphism. On $\nu_pM$ we have the isoparametric foliation $\F_p$ and on $T_pM$ the trivial foliation by $\{T_pM\}$. By construction $\phi$ maps leaves of the product foliation on $T_pM\times\nu_p M$ diffeomorphically onto leaves of $\F$. 
\end{proof}
\begin{remark}
Following the proof we can easily show the above product is metric if $X$ is Euclidean. Let $X$ be the Euclidean space. Then a leaf $M$ of lowest dimension is an affine subspace, because its normal exponential map is a diffeomorphism, and its global slices are the orthogonal, complementary affine subspaces. Let $\F'$ denote the foliation by parallel affine subspaces of $M$. As stated in the proof each other leaf is a union of leaves of $\F'$. This means the orthogonal projection of $X$ onto a given slice $S$ of $M$ preserves the leaves of $\F'$ and therefore the leaves of $\F$. Together with the orthogonal projection onto $M$ we have a {\it metric} foliated decomposition $X=M\times S$. This recovers the result about isoparametric splitting; see Corollary 6.3.12. in \cite{palaisterng}. 
\end{remark}


\end{document}